\documentclass[oneside,english]{amsart}
\usepackage[T1]{fontenc}
\usepackage[latin9]{inputenc}
\usepackage{textcomp}
\usepackage{amsbsy}
\usepackage{amstext}
\usepackage{amsthm}
\usepackage{amssymb}

\makeatletter
\numberwithin{equation}{section}
\numberwithin{figure}{section}
\theoremstyle{plain}
\newtheorem*{thm*}{\protect\theoremname}
\theoremstyle{plain}
\newtheorem{thm}{\protect\theoremname}
\theoremstyle{definition}
\newtheorem{defn}[thm]{\protect\definitionname}
\theoremstyle{remark}
\newtheorem{rem}[thm]{\protect\remarkname}
\theoremstyle{plain}
\newtheorem{prop}[thm]{\protect\propositionname}
\theoremstyle{plain}
\newtheorem{lem}[thm]{\protect\lemmaname}
\theoremstyle{plain}
\newtheorem{cor}[thm]{\protect\corollaryname}

\usepackage[a4paper]{geometry}

\theoremstyle{definition}

\theoremstyle{problem}

\DeclareMathOperator{\II}{\mathrm{I\!I}}

\DeclareMathOperator{\tr}{tr}

\DeclareMathOperator{\GL}{GL}

\DeclareMathOperator{\Fr}{Fr}

\DeclareMathOperator{\Ric}{Ric}

\DeclareMathOperator{\HH}{H}

\usepackage{graphicx}
\usepackage{amsmath}
\usepackage{tikz}
\usepackage{tikz-3dplot}
\usetikzlibrary{math,calc,arrows}

\makeatother

\usepackage{babel}
\providecommand{\corollaryname}{Corollary}
\providecommand{\definitionname}{Definition}
\providecommand{\lemmaname}{Lemma}
\providecommand{\propositionname}{Proposition}
\providecommand{\remarkname}{Remark}
\providecommand{\theoremname}{Theorem}

\begin{document}
\title[Isometric Embeddings of AH Manifolds into Hyperbolic Space]{Isometric Embeddings of Asymptotically Hyperbolic Manifolds into
Hyperbolic Space}
\author{Marco Usula}
\begin{abstract}
The Nash Embedding theorem asserts that every closed Riemannian manifold
can be isometrically embedded into a Euclidean space of sufficiently
high dimension. In this paper, continuing the analysis initiated in
\cite{UsulaIsometricEmbeddings}, we address the asymptotically hyperbolic
version of the problem. We prove that an asymptotically hyperbolic
manifold $\left(M,g\right)$ can be realized as a submanifold of a
sufficiently high-dimensional hyperbolic space $\HH^{N+1}$, intersecting
the sphere at infinity transversely, if and only if $g$ is asymptotically
Poincaré--Einstein to second order.
\end{abstract}

\maketitle

\section{Introduction}

The Nash Embedding Theorem is a foundational result in Riemannian
geometry. The theorem asserts that every closed Riemannian manifold
can be isometrically embedded into a Euclidean space of sufficiently
high dimension. This result was proved by John Nash across the years
1954--1956 in two Annals papers \cite{NashC1IsometricImbeddings, NashImbeddingProblem}:
the former proves the existence of $C^{1}$ embeddings, while the
latter proves the existence of $C^{\infty}$ embeddings. Besides the
intrinsic importance of the Nash Embedding Theorem, \cite{NashImbeddingProblem}
also introduces a fundamental analytic tool, the \emph{Nash--Moser
Technique}, which has been used to prove existence theorems for several
non-linear partial differential equations. We refer to \cite{GuntherIsometricEmbeddings, GuntherPerturbationProblem, YangNashIsometricEmbedding, TaoNashEmbedding, GromovPDR, GromovDescendants}
for extensive discussions on the Nash Embedding Theorem, the Nash--Moser
Technique, and the vast amount of mathematical literature stemming
from them.

In this paper, continuing the work initiated by the author in \cite{UsulaIsometricEmbeddings},
we address a similar isometric embedding problem in the context of
\emph{conformally compact geometry}. Given a compact manifold with
boundary $M^{m}$, a \emph{conformally compact metric} on $M$ is
a smooth metric $g$ on the interior $M^{\circ}$ such that, given
a boundary defining function\footnote{A boundary defining function for $M$ is a smooth function $\rho:M\to[0,+\infty)$
such that $\rho^{-1}\left(0\right)=\partial M$ and $0$ is a regular
value of $\rho$.} $\rho$ for $M$, the conformal rescaling $\rho^{2}g$ extends smoothly
to a metric on $M$. This definition mimics the behavior of the metric
of hyperbolic space. Indeed, in its Poincaré model, the hyperbolic
$m$-dimensional space $\HH^{m}$ is the closed unit ball in $\mathbb{R}^{m}$
with coordinates $\boldsymbol{y}=\left(y^{1},...,y^{m}\right)$, equipped
with the conformally compact metric
\[
g_{\mathrm{hyp}}=\frac{4\left|d\boldsymbol{y}\right|^{2}}{\left(1-\left|\boldsymbol{y}\right|^{2}\right)^{2}}.
\]
The study of the geometry and analysis on conformally compact geometry
was pioneered by Mazzeo, Melrose, Fefferman, Graham, Lee, Biquard
and many others; we refer to the foundational papers \cite{MazzeoPhD, MazzeoHodge, MazzeoMelroseResolvent, FeffermanGraham, GrahamLee, LeeFredholm, Biquard_Metriques}
for many more details on this subject.

The most important class of conformally compact metrics is arguably
that of \emph{asymptotically hyperbolic} metrics. As the name suggests,
a conformally compact metric is asymptotically hyperbolic if its sectional
curvature at infinity $\kappa_{\infty}^{g}$ (cf. Definition \ref{def:sectional-curvature-at-infinity})
is constant and equal to $-1$. This class of manifolds has been studied
extensively, particularly in the context of the Einstein equation.
If $\left(M^{m},g\right)$ is conformally compact and Einstein, then
it is forced to have constant sectional curvature at infinity; the
common normalization is to require this constant to be $-1$, so that
$g$ must satisfy the non-linear partial differential equation
\[
\Ric g+\left(m-1\right)g=0.
\]
The conformally compact metrics satisfying this equation are called
\emph{Poincaré--Einstein }metrics, and they are central in the study
of invariants of closed \emph{conformal} manifolds, and in the study
of the Riemannian version of the AdS-CFT correspondence \cite{FeffermanGraham, Maldacena, WittenHolography, GrahamLee, LeeFredholm, BiquardAdS, UsulaYM, LimaEinsteinYangMills}
.

From a conformally compact perspective, it is natural to consider
compact submanifolds with boundary of a conformally compact manifold
$\left(M,g\right)$ intersecting the boundary at infinity transversely
(the so-called \emph{p-submanifolds}, cf. Definition \ref{def:p-submanifolds}).
Indeed, if $S\subset M$ is such a submanifold, then the induced metric
on $S$ is conformally compact (cf. Proposition 3.6 of \cite{UsulaIsometricEmbeddings}).
This leads to the following natural question: when can we realize
a conformally compact manifold as a p-submanifold of a sufficiently
high-dimensional hyperbolic space? The author partially solved this
problem in \cite{UsulaIsometricEmbeddings}; there, it is explained
that the answer depends fundamentally on the sectional curvature at
infinity $\kappa_{\infty}^{g}$ of $g$. Specifically, it is proved
that if $\left(M^{m},g\right)$ can be isometrically p-embedded in
$\mathrm{H}^{N+1}$ for some $N\geq m$, then $\kappa_{\infty}^{g}$
must satisfy the pointwise inequality
\[
\kappa_{\infty}^{g}\geq-1.
\]
Moreover, the main result of \cite{UsulaIsometricEmbeddings} is that
if the stronger inequality
\[
\kappa_{\infty}^{g}>-1
\]
holds, then the condition is also sufficient.

\cite{UsulaIsometricEmbeddings} leaves open the following natural
question concerning the asymptotically hyperbolic case: given an asymptotically
hyperbolic manifold $\left(M,g\right)$ (that is, $\kappa_{\infty}^{g}\equiv-1$),
can we p-embed it isometrically into $\mathrm{H}^{N+1}$ for some
$N$ large enough? This paper answers this question definitively:
somewhat surprisingly, this is possible if and only if $g$ is \emph{asymptotically
Poincaré--Einstein }to second order. More precisely, the main result
of the paper is the following
\begin{thm*}
Let $\left(M^{m},g\right)$ be an asymptotically hyperbolic manifold,
with $m\geq2$. Then the following are equivalent:
\begin{enumerate}
\item there exists $N\in\mathbb{N}$, $N\geq m$, and an embedding $u:M\to\mathrm{H}^{N+1}$
mapping boundary to boundary, interior to interior, and transverse
to $\partial\mathrm{H}^{N+1}$, such that $u^{*}g_{\mathrm{hyp}}=g$;
\item $g$ is asymptotically Poincaré--Einstein to second order: $\left|\Ric g+\left(m-1\right)g\right|_{g}=O\left(x^{2}\right)$,
where $x$ is any boundary defining function for $M$.
\end{enumerate}
Moreover, if the condition $\left|\Ric g+\left(m-1\right)g\right|_{g}=O\left(x^{2}\right)$
is satisfied, then $N$ can be taken to be $\mathrm{N}_{m+1}$, where
$\mathrm{N}_{k}$ is the smallest positive integer such that every
closed Riemannian $k$-manifold can be isometrically embedded in $\mathbb{R}^{\mathrm{N}_{k}}$.
\end{thm*}
The proof of this result uses fundamentally the Nash Embedding Theorem,
and some of the ideas presented in \cite{UsulaIsometricEmbeddings}.
More precisely, the idea of the main proof in \cite{UsulaIsometricEmbeddings}
is the following: given $\left(M^{m},g\right)$ conformally compact
and with $\kappa_{\infty}^{g}>-1$, one can construct a boundary defining
function $x$ for $M$ such that the symmetric $2$-tensor $G=x^{2}g-dx^{2}$
is a smooth metric on $M$. Given such a $x$, the Nash Embedding
Theorem implies that there exists an isometric embedding $\boldsymbol{v}:\left(M,G\right)\to\mathbb{R}^{N}$
for some $N\geq m$. The map $\left(x,\boldsymbol{v}\right):M\to[0,+\infty)\times\mathbb{R}^{N}$
is then an isometric p-embedding from $\left(M,g\right)$ to the half-space
model of $\HH^{N+1}$.

When $\kappa_{\infty}^{g}\equiv-1$, a naive implementation of the
same strategy fails: indeed, for \emph{any }boundary defining function
$x$, the symmetric $2$-tensor $G=x^{2}g-dx^{2}$ on $M$ is necessarily
degenerate along $\partial M$ in the $\partial_{x}$ direction. \emph{A
posteriori}, one cannot expect every asymptotically hyperbolic $\left(M,g\right)$
to admit an isometric p-embedding into some hyperbolic space: indeed,
if $M$ is a p-submanifold of $\HH^{N+1}$, then it is easy to see
that $M$ is asymptotically hyperbolic if and only if it intersects
the sphere at infinity of $\HH^{N+1}$ \emph{orthogonally}. This immediately
implies that the conformal infinity of $\partial M$ must be \emph{totally
umbilic} in $M$ (that is, the trace-free part of the second fundamental
form of $\partial M$ in $M$ must vanish). With some extra work (cf.
Proposition \ref{prop:necessity}), we show that something stronger
is true: $\partial M$ must be \emph{totally geodesic} in $M$ in
any\emph{ geodesic compactification} (cf. Definition \ref{def:geodesic-compactification})
of $\left(M,g\right)$. This condition is equivalent to the fact that
$\left|\Ric g+\left(m-1\right)g\right|_{g}$ vanishes to second order
along $\partial M$, namely that $g$ is asymptotically Poincaré--Einstein
to second order. Somewhat surprisingly, this condition is also \emph{sufficient},
as we show in §\ref{subsec:Sufficiency}. The idea of the proof is
that, even if for every boundary defining function $x$ the symmetric
$2$-tensor $G=x^{2}g-dx^{2}$ is degenerate in the $\partial_{x}$
direction, thanks to the condition $\left|\Ric g+\left(m-1\right)g\right|_{g}=O\left(x^{2}\right)$
we can construct $x$ in such a way that:
\begin{enumerate}
\item $G$ is positive-definite \emph{in the interior}, and
\item $G$ can be realized as the pull-back of a genuine smooth metric $\mathcal{G}$
on the manifold $\mathcal{M}=M\times S^{1}$, via a map $F:M\to\mathcal{M}$
which restricts to an embedding of the interiors and the boundaries,
but whose differential is degenerate in the direction orthogonal to
$\partial M$ in $M$.
\end{enumerate}
The isometric p-embedding of $\left(M,g\right)$ into $\HH^{N+1}$
is then realized as $\left(x,\boldsymbol{w}\right)$, where $\boldsymbol{w}=\boldsymbol{v}\circ F$
and $\boldsymbol{v}:\left(\mathcal{M},\mathcal{G}\right)\to\mathbb{R}^{N}$
is an isometric embedding, whose existence is guaranteed by the Nash
Embedding Theorem.

\subsection*{Acknowledgements}

The author is grateful to the Wiener--Anspach Foundation for the
support.

\section{\label{sec:Conformally-compact-geometry}Conformally compact geometry}

In this section, we recall the basic notions of conformally compact
geometry, which we shall use in the rest of the paper. For an extensive
discussion on these topics, we refer to \cite{FeffermanGraham, BiquardAdS, MazzeoPhD, GrahamLee}.

We shall work on a compact manifold with boundary $M$; we denote
by $\partial M$ its boundary, and by $M^{\circ}$ its interior. Recall
that a \emph{boundary defining function }for $M$ is a smooth function
$x:M\to[0,+\infty)$ such that $x^{-1}\left(0\right)=\partial M$
and $dx$ is nowhere zero along $\partial M$.
\begin{defn}
A \emph{conformally compact metric} on $M$ is a metric $g$ on $M^{\circ}$
such that, for some (hence every) boundary defining function $\rho$
on $M$, the conformal rescaling $\rho^{2}g$ extends smoothly from
$M^{\circ}$ to a metric on $M$.
\end{defn}

We emphasize the fact that a conformally compact metric on $M$ is
really a smooth object on the whole of $M$, not just on the interior
$M^{\circ}$. To formalize this concept, we introduce a central object
in conformally compact geometry: the \emph{$0$-tangent bundle} $^{0}TM$.
This bundle is obtained, abstractly via the Serre--Swan Theorem,
as the bundle whose $C^{\infty}\left(M\right)$ module of sections
is the space $\mathcal{V}_{0}\left(M\right)$ of vector fields on
$M$ vanishing along $\partial M$. The inclusion $\mathcal{V}_{0}\left(M\right)\to\mathcal{V}\left(M\right)$
induces a smooth bundle map $\#:{}^{0}TM\to TM$, and since $0$-vector
fields are unconstrained in the interior, $\#$ restricts to an isomorphism
$^{0}TM_{|M^{\circ}}\to TM^{\circ}$. However, due to the uniform
degeneracy of $0$-vector fields along $\partial M$, $\#$ restricts
to the zero map along the boundary.

If $E$ is a tensor bundle for $M$ associated to the principal $\GL\left(m,\mathbb{R}\right)$
bundle $\Fr\left(TM\right)$ via a real linear representation of $\GL\left(m,\mathbb{R}\right)$,
we can construct its ``$0$-version'' $^{0}E$, associated to $\Fr\left(^{0}TM\right)$
via the same representation. The canonical isomorphism ${^{0}TM}_{|M^{\circ}}\simeq TM^{\circ}$
then induces a canonical isomorphism $^{0}E_{M^{\circ}}\simeq E_{|M^{\circ}}$,
so that we can think of sections of $^{0}E$ as sections of $E_{|M^{\circ}}$
exhibiting a uniform degeneracy or singularity along $\partial M$,
according to the degrees of covariance / contravariance of $E$. More
concretely, in this paper we shall use the following associated $0$-tensor
bundles:
\begin{enumerate}
\item the \emph{$0$-cotangent bundle }$^{0}T^{*}M$, namely the dual of
$^{0}TM$, and its $k$-th exterior power $^{0}\Lambda^{k}$ called
the bundle of \emph{$0$-$k$-forms}; smooth sections of $^{0}\Lambda^{k}$
are precisely $k$-forms $\omega$ on $M^{\circ}$ such that $x^{k}\omega$
extends smoothly to a $k$-form on $M$;
\item the bundle $S^{2}\left(^{0}T^{*}M\right)$ of \emph{symmetric $0$-$2$-tensors};
smooth sections of $S^{2}\left(^{0}T^{*}M\right)$ are symmetric $2$-tensors
$g$ on $M^{\circ}$ such that $x^{2}g$ extends smoothly to a symmetric
$2$-tensor on $M$.
\end{enumerate}
From the previous point, we see immediately that conformally compact
metrics are precisely fibrewise positive-definite sections of $S^{2}\left(^{0}T^{*}M\right)$,
that is, bundle metrics on $^{0}TM$. For this reason, conformally
compact metrics are often called \emph{$0$-metrics}. Note that a
$0$-metric on $M$ induces a bundle metric on any $0$-tensor bundle
$^{0}E$; therefore, by construction, if $s$ is a smooth section
of $^{0}E$ then its pointwise norm $\left|s\right|_{g}^{2}$ extends
from the interior to a smooth function on $M$.

Conformally compact manifolds share many properties with hyperbolic
space. First of all, if $\left(M,g\right)$ is conformally compact,
then the metric on $M^{\circ}$ induced by $g$ is complete, with
bounded geometry and infinite volume. Moreover, while $g$ does not
extend to a metric over $\partial M$, it does induce a \emph{conformal
class}, defined as
\[
\mathfrak{c}_{\infty}\left(g\right)=\left\{ \left(x^{2}g\right)_{|\partial M}:\text{\ensuremath{x} is a boundary defining function on \ensuremath{M}}\right\} .
\]
This is called the \emph{conformal infinity }of $g$. Finally, $g$
is \emph{asymptotically negatively curved}. To explain this more precisely,
we introduce the following invariant of $g$ at infinity:
\begin{defn}
\label{def:sectional-curvature-at-infinity}Let $x$ be a boundary
defining function for $M$. The \emph{sectional curvature at infinity
}of $g$ is the function
\[
\kappa_{\infty}^{g}:=-\left.\left|\frac{dx}{x}\right|_{g}^{2}\right|_{M}\in C^{\infty}\left(\partial M\right).
\]
\end{defn}

\begin{rem}
The definition above is well-posed, meaning that it does not depend
on $x$. Indeed, if $x'=e^{\varphi}x$ is another boundary defining
function for $M$, then
\[
\frac{dx'}{x'}=\frac{dx}{x}+d\varphi.
\]
Now $d\varphi$ is a smooth section of $\Lambda^{1}$, and therefore
it vanishes identically along $\partial M$ when interpreted as a
smooth section of $^{0}\Lambda^{1}$ vanishing along $\partial M$.
\end{rem}

The name of $\kappa_{\infty}^{g}$ is due to the following property:
\begin{prop}
(Proposition 1.13 of \cite{MazzeoPhD}) If $\left\{ p_{k}\right\} _{k\in\mathbb{N}}$
is a sequence of points in $M^{\circ}$ converging as $k\to\infty$
to a point $p_{\infty}\in\partial M$, and $\pi_{k}\subset T_{p_{k}}M^{\circ}$
is a sequence of tangent planes, then the sectional curvature $\kappa^{g}\left(\pi_{k}\right)$
converges to $\kappa_{\infty}^{g}\left(p_{\infty}\right)$ as $k\to\infty$.
\end{prop}

This property justifies the following
\begin{defn}
An \emph{asymptotically hyperbolic }metric on $M$ is a conformally
compact metric $g$ such that $\kappa_{\infty}^{g}\equiv-1$.
\end{defn}

The class of asymptotically hyperbolic and \emph{Einstein }metrics
is particularly well-studied. These metrics are called \emph{Poincaré--Einstein}.
If $\left(M^{m},g\right)$ is asymptotically hyperbolic and Einstein,
the condition $\kappa_{\infty}^{g}\equiv-1$ implies that the scalar
curvature of $g$ must equal that of $\HH^{m}$. For this reason,
$g$ is Poincaré--Einstein if and only if it satisfies the non-linear
PDE
\[
\Ric g+\left(m-1\right)g=0.
\]
Poincaré--Einstein metrics are central in the study of the Riemannian
version of the AdS-CFT correspondence (cf. \cite{BiquardAdS}), and
in the study of conformal invariants of closed manifolds (cf. \cite{FeffermanGraham}).

\section{\label{sec:p-submanifolds-of-hyperbolic}p-submanifolds of hyperbolic
spaces}
\begin{defn}
\label{def:p-submanifolds}Let $S,M$ be compact manifolds with boundary.
An \emph{interior p-embedding} \emph{of $S$ in $M$} is an injective
immersion $\iota:S\to M$ such that, for some (hence every) boundary
defining function $\rho$ for $M$, $\iota^{*}\rho$ is a boundary
defining function for $S$. An \emph{interior p-submanifold} of $M$
is the image of an interior p-embedding.
\end{defn}

\begin{rem}
Note that, if $\iota:S\to M$ is an interior p-embedding, then $\iota$
maps $\partial S$ to $\partial M$ and $S^{\circ}$ to $M^{\circ}$,
and both restrictions are embeddings. We can equivalently define an
interior p-submanifold as a topologically closed submanifold $S\subset M$
which intersects $\partial M$ transversely.
\end{rem}

The class of interior p-submanifolds is arguably the most natural
class of submanifolds of a conformally compact manifold to consider,
because of the following
\begin{prop}
(Proposition 3.6 in \cite{UsulaIsometricEmbeddings}) Let $\left(M,g\right)$
be a conformally compact manifold, and let $\iota:S\to M$ be an interior
p-embedding. Then $\iota^{*}g$ is a conformally compact metric on
$S$. Moreover, the sectional curvatures at infinity of $g$ and $\iota^{*}g$
satisfy the inequality
\[
\kappa_{\infty}^{\iota^{*}g}\geq\kappa_{\infty|\partial S}^{g}.
\]
\end{prop}

Given $\lambda^{2}\in\mathbb{R}^{+}$, denote by $\HH^{N+1}\left(-\lambda^{2}\right)$
the rescaled hyperbolic space with constant sectional curvature $-\lambda^{2}$
(when $\lambda^{2}=1$, we simply write $\HH^{N+1}$ in place of $\HH^{N+1}\left(-1\right)$).
In \cite{UsulaIsometricEmbeddings}, the author considered the following
isometric embedding problem. Given a conformally compact manifold
$\left(M^{m},g\right)$, can we realize it as an interior p-submanifold
of $\HH^{N+1}\left(-\lambda^{2}\right)$? The previous proposition
immediately shows that a necessary condition is the inequality
\[
\kappa_{\infty}^{g}\geq-\lambda^{2}.
\]
The main result of \cite{UsulaIsometricEmbeddings} is that, if the
previous inequality is \emph{strict}, then the condition is also sufficient:
\begin{thm}
\label{thm:previous-main-theorem}(Theorem 4.2 of \cite{UsulaIsometricEmbeddings})
Let $\left(M^{m},g\right)$ be a conformally compact manifold, and
let $\lambda^{2}\in\mathbb{R}^{+}$. Assume that $\kappa_{\infty}^{g}>-\lambda^{2}$.
Then $\left(M,g\right)$ can be isometrically p-embedded into $\HH^{\mathrm{N}_{m}+1}\left(-\lambda^{2}\right)$,
where $\mathrm{N}_{k}$ is the smallest positive integer such that
every closed Riemannian $k$-manifold can be isometrically embedded
in $\mathbb{R}^{\mathrm{N}_{k}}$.
\end{thm}

\section{The main theorem}

The natural question left open in \cite{UsulaIsometricEmbeddings}
is: if $\left(M,g\right)$ only satisfies the non-strict inequality
$\kappa_{\infty}^{g}\geq-\lambda^{2}$, is it possible to p-embed
$\left(M,g\right)$ isometrically into $\HH^{N+1}\left(-\lambda^{2}\right)$
for $N$ sufficiently large? In this paper, we solve this problem
definitively in the most relevant special case, namely that of \emph{asymptotically
hyperbolic} submanifolds of hyperbolic space. If $\left(M,g\right)\subset\HH^{N+1}$
is an interior p-submanifold, then $\kappa_{\infty}^{g}\geq-1$, with
equality if and only if $M$ intersects the sphere at infinity \emph{orthogonally
}(a proof of this fact will be given below). The main result of this
paper is that, somewhat surprisingly, $\left(M,g\right)$ can be isometrically
p-embedded into some hyperbolic space $\HH^{N+1}$ if and only if
$g$ is \emph{asymptotically Poincaré--Einstein to second order}.
The precise statement is the following:
\begin{thm}
\label{thm:main-thm}Let $\left(M^{m},g\right)$ be an asymptotically
hyperbolic manifold, with $m\geq2$. Then $\left(M,g\right)$ can
be isometrically p-embedded into $\HH^{N+1}$ for some $N\geq m$
if and only if the symmetric $0$-$2$-tensor $\Ric g+\left(m-1\right)g$
vanishes to second order along $\partial M$, namely
\[
\Ric g+\left(m-1\right)g=O\left(x^{2}\right)
\]
where $x$ is a boundary defining function for $M$. In this case,
$\left(M,g\right)$ can be isometrically p-embedded into $\HH^{\mathrm{N}_{m+1}+1}$,
where $\mathrm{N}_{k}$ is the smallest positive integer such that
every closed Riemannian $k$-manifold embeds isometrically into $\mathbb{R}^{\mathrm{N}_{k}}$.
\end{thm}

The rest of the paper is devoted to the proof of this result.

\subsection{\label{subsec:Necessity}Necessity}

In this sub-section, we prove the necessity part of Theorem \ref{thm:main-thm};
namely, we prove the following
\begin{prop}
\label{prop:necessity}Let $M^{m}$ be an interior p-submanifold of
$\HH^{N+1}$, with $m\geq2$ and $N\geq m$. If the induced conformally
compact metric $g$ on $M$ is asymptotically hyperbolic, then
\[
\Ric g+\left(m-1\right)g=O\left(x^{2}\right)
\]
as a section of $S^{2}\left(^{0}T^{*}M\right)$, where $x$ is a boundary
defining function for $M$.
\end{prop}

The proof requires some more known facts about conformally compact
geometry. For a generic conformally compact metric $g$ on $M$, the
symmetric $0$-$2$-tensor $\Ric g+\left(m-1\right)g$ can be non-zero
along the boundary. As observed by Graham and Lee (cf. Lemma 2.1 in
\cite{GrahamLee}), the vanishing condition $\Ric g+\left(m-1\right)g=O\left(x\right)$
is equivalent to the condition that $g$ is asymptotically hyperbolic.
Graham and Lee then proved (cf. Lemma 5.2 in \cite{GrahamLee}) that,
if $g$ is asymptotically hyperbolic, then for every metric $h_{0}$
in the conformal infinity $\mathfrak{c}_{\infty}\left(g\right)$ of
$g$ we can construct a so-called \emph{geodesic} \emph{boundary defining
function} $r$, with the following properties:
\begin{enumerate}
\item $\left(r^{2}g\right)_{|\partial M}=h_{0}$;
\item $\left|dr\right|_{r^{2}g}\equiv1$ on a neighborhood of $\partial M$.
\end{enumerate}
Such a boundary defining function $r$ is unique in a neighborhood
of $\partial M$. Using the flow of the gradient of $r$ with respect
to the compactified metric $r^{2}g$, we can construct a\emph{ }collar
$[0,\varepsilon)\times\partial M\to M$ where $\left|dr\right|_{r^{2}g}\equiv1$;
in this collar, $g$ takes the normal form
\[
g=\frac{dr^{2}+h\left(r\right)}{r^{2}},
\]
where $h:[0,\varepsilon)\to C^{\infty}\left(\partial M;S^{2}\left(T^{*}\partial M\right)\right)$
is a smooth family of metrics on $\partial M$ such that $h\left(0\right)=h_{0}$.
\begin{defn}
\label{def:geodesic-compactification}Let $\left(M,g\right)$ be asymptotically
hyperbolic. A \emph{geodesic compactification} of $g$ is a conformal
compactification $\overline{g}=r^{2}g$, where $r$ is a geodesic
boundary defining function for $\left(M,g\right)$.
\end{defn}

Consider now the asymptotic expansion of $h$ in powers of $r$:
\[
h\sim\sum_{j=0}^{\infty}h_{j}r^{j}.
\]
The first coefficient $h_{1}=\dot{h}\left(0\right)$ is related to
the second fundamental form $\II$ of $\left(\partial M,h_{0}\right)$
in $\left(M,r^{2}g\right)$ by the formula
\[
h_{1}=-2\II;
\]
A proof of this fact is given for example in Lemma A.16 of \cite{UsulaNahmPoles}.
Finally, we recall Formula (4) of \cite{MazzeoPacardMaskit}\footnote{In \cite{MazzeoPacardMaskit} this formula is expressed by saying
that $\Ric g+\left(m-1\right)g$ equals the singular term above plus
a $O\left(1\right)$ section of $S^{2}\left(T^{*}M\right)$; this
is equivalent to our formulation, because smooth sections of $S^{2}\left(T^{*}M\right)$
are precisely $O\left(r^{2}\right)$ sections of $S^{2}\left(^{0}T^{*}M\right)$,
cf. §2.}
\[
\Ric g+\left(m-1\right)g=\frac{1}{2r}\left(\left(m-2\right)h_{1}+\left(\tr_{h_{0}}h_{1}\right)\left(dr^{2}+h_{0}\right)\right)+O\left(r^{2}\right).
\]
From this formula, it follows that $\Ric g+\left(m-1\right)g=O\left(r^{2}\right)$
if and only if 
\begin{align*}
\left(m-2\right)h_{1} & =0\\
\tr_{h_{0}}h_{1} & =0.
\end{align*}
If $m>2$, the first condition implies $h_{1}=0$. If $m=2$, the
first equation becomes trivial, but since $h_{1}$ is pure trace because
$\partial M$ is $1$-dimensional, the condition $\tr_{h_{0}}h_{1}=0$
is equivalent to $h_{1}=0$. In both cases, since $h_{1}=-2\II$,
we have that $\Ric g+\left(m-1\right)g=O\left(r^{2}\right)$ if and
only if $\left(\partial M,h_{0}\right)$ is totally geodesic in $\left(M,r^{2}g\right)$.

We are almost ready to prove Proposition \ref{prop:necessity}. Let
$M^{m}$ be an interior p-submanifold of $\HH^{N+1}$, with $m\geq2$
and $N\geq m$. Since $M$ is compact, $\partial M$ is a compact
submanifold of the sphere at infinity of $\HH^{N+1}$; since $N\geq m$,
$\partial M$ has codimension at least $1$ in $S^{N}=\partial\HH^{N+1}$.
Therefore, $\partial M$ misses at least a point of $S^{N}$. For
this reason, without loss of generality, we can assume that $M$ is
p-embedded in the upper half-space model of $\HH^{N+1}$. More precisely,
we can think of $\HH^{N+1}$ as the topologically closed manifold
with boundary $\mathbb{R}_{1}^{N+1}:=[0,+\infty)\times\mathbb{R}^{N}$,
with coordinates $\left(X,\boldsymbol{Y}\right)=\left(X,Y^{1},...,Y^{N}\right)$
and with hyperbolic $0$-metric
\[
g_{\mathrm{hyp}}=\frac{dX^{2}+\left|d\boldsymbol{Y}\right|^{2}}{X^{2}}.
\]
Denote by $x$ and $\boldsymbol{y}$ the pull-backs of $X$ and $\boldsymbol{Y}$
to $M$. Then the conformally compact metric $g$ on $M$ induced
by the interior p-embedding $\iota:M\to\HH^{N+1}$ takes the form
\[
g=\frac{dx^{2}+\left|d\boldsymbol{y}\right|^{2}}{x^{2}}.
\]
Since $M$ is an interior p-submanifold of $\mathbb{R}_{1}^{N+1}$,
$x$ is a boundary defining function for $M$. Note also that the
symmetric $2$-tensor $\left|d\boldsymbol{y}\right|^{2}$ can contain
terms involving the $1$-form $dx$, so $g$ is not necessarily asymptotically
hyperbolic.
\begin{lem}
\label{lem:x-is-asymptotically-geodesic-bdf}The conformally compact
metric $g$ is asymptotically hyperbolic if and only if $M$ intersects
the boundary at infinity of $\HH^{N+1}$ orthogonally. In this case,
we have
\[
\left|dx\right|_{x^{2}g}^{2}=1+O\left(x^{2}\right).
\]
\end{lem}

\begin{rem}
Even if the metric on $\HH^{N+1}$ is not well-defined up to the boundary,
its conformal structure is. Therefore, the notion of orthogonality
at infinity is well-defined.
\end{rem}

\begin{proof}
Define $\overline{g}=x^{2}g=dx^{2}+\left|d\boldsymbol{y}\right|^{2}$,
and denote by $\overline{\nabla}x$ the gradient of $x$ with respect
to $\overline{g}$. Consider now the vector field $\partial_{X}$
on $\HH^{N+1}$; the restriction of $\partial_{X}$ to $M$ can be
decomposed as
\[
\partial_{X}=\overline{\nabla}x+\left(\partial_{X}\right)^{\bot},
\]
where $\left(\partial_{X}\right)^{\bot}$ is the section of $T\left(\mathbb{R}_{1}^{N+1}\right)_{|M}$
obtained by projecting $\partial_{X}$ to the orthogonal complement
of $TM$ with respect to the Euclidean metric on $\mathbb{R}_{1}^{N+1}$.
We then have
\[
\left|dx\right|_{\overline{g}}^{2}=\left|\overline{\nabla}x\right|_{\overline{g}}^{2}=1-\left|\left(\partial_{X}\right)^{\bot}\right|_{\mathbb{R}_{1}^{N+1}}^{2}.
\]
Now, $M$ intersects the boundary of $\mathbb{R}_{1}^{N+1}$ orthogonally
if and only if $\partial_{X}$ is tangent to $M$ along $\partial M$,
or equivalently if and only if $\left(\partial_{X}\right)^{\bot}\equiv0$
along $\partial M$. In this case, we have $\left|\left(\partial_{X}\right)^{\bot}\right|_{\mathbb{R}_{1}^{N+1}}^{2}=O\left(x^{2}\right)$
and therefore $\left|dx\right|_{\overline{g}}^{2}=1+O\left(x^{2}\right)$
as claimed. 
\end{proof}
\begin{cor}
\label{cor:equivalence-of-second-fundamental-forms}Let $r$ be a
geodesic boundary defining function for $\left(M,g\right)$ inducing
the metric $h_{0}=\left(x^{2}g\right)_{|\partial M}$. Then the second
fundamental forms of $\left(\partial M,h_{0}\right)$ in the compactified
metrics $\left(M,x^{2}g\right)$ and $\left(M,r^{2}g\right)$ coincide.
\end{cor}

\begin{proof}
Write $x=e^{\varphi}r$ for some $\varphi\in C^{\infty}\left(M\right)$.
Since $x$ and $r$ induce the same metric in the conformal infinity
of $g$, we have 
\[
\varphi_{|\partial M}\equiv0.
\]
Now, we have
\begin{align*}
\left|dx\right|_{x^{2}g}^{2} & =\left|\frac{dx}{x}\right|_{g}^{2}=\left|\frac{dr}{r}\right|_{g}^{2}+\left|d\varphi\right|_{g}^{2}+2\left\langle \frac{dr}{r},d\varphi\right\rangle _{g}.
\end{align*}
Since $d\varphi$ is a smooth $1$-form on $M$, it is $O\left(r\right)$
as a section of $^{0}T^{*}M$, and therefore $\left|d\varphi\right|_{g}^{2}=O\left(r^{2}\right)$.
Since $r$ is a geodesic boundary defining function for $\left(M,g\right)$,
we have $\left|dr/r\right|_{g}^{2}=\left|dr\right|_{r^{2}g}^{2}\equiv1$
on a geodesic collar induced by $r$. Furthermore, in such a geodesic
collar, we have
\begin{align*}
\left\langle \frac{dr}{r},d\varphi\right\rangle _{g} & =r\left\langle dr,d\varphi\right\rangle _{r^{2}g}\\
 & =r\partial_{r}\varphi.
\end{align*}
Substituting these conditions and the condition $\left|dx\right|_{x^{2}g}^{2}=1+O\left(r^{2}\right)$
in the previous identity, we obtain $r\partial_{r}\varphi=O\left(r^{2}\right)$,
which is equivalent to 
\[
\left(\partial_{r}\varphi\right)_{|\partial M}\equiv0.
\]
Finally, since $x^{2}g=e^{2\varphi}\left(r^{2}g\right)$, the second
fundamental forms $\II_{x^{2}g}$ and $\II_{r^{2}g}$ of $\partial M$
in $\left(M,x^{2}g\right)$ and $\left(M,r^{2}g\right)$ are related
by the equation
\[
\II_{x^{2}g}=e^{\varphi_{|\partial M}}\left(\II_{r^{2}g}-\left(\partial_{r}\varphi\right)_{|\partial M}h_{0}\right).
\]
Substituting $\varphi_{|\partial M}\equiv0$ and $\left(\partial_{r}\varphi\right)_{|\partial M}\equiv0$,
we obtain the claimed identity $\II_{x^{2}g}=\II_{r^{2}g}$.
\end{proof}
We can now complete the proof of Proposition \ref{prop:necessity}:
\begin{proof}[Proof of Proposition \ref{prop:necessity}]
As in the proof of Lemma \ref{lem:x-is-asymptotically-geodesic-bdf},
define $\overline{g}=x^{2}g=dx^{2}+\left|d\boldsymbol{y}\right|^{2}$,
so that $\left(M,\overline{g}\right)$ is an interior p-submanifold
of $\mathbb{R}_{1}^{N+1}$ meeting the boundary orthogonally. By Corollary
\ref{cor:equivalence-of-second-fundamental-forms} and the discussion
at the beginning of this sub-section, it suffices to prove that $\partial M$
is totally geodesic in $\left(M,\overline{g}\right)$. Denote by $\nabla$
the Levi-Civita connection on $\mathbb{R}_{1}^{N+1}$, and denote
by $\overline{\nabla}$ the Levi-Civita connection on $\left(M,\overline{g}\right)$.
$\partial M$ is totally geodesic in $\left(M,\overline{g}\right)$
if and only if, for every pair of vector fields $U,V$ on $M$ tangent
to $\partial M$, $\overline{\nabla}_{U}V$ is still tangent to $\partial M$.
Now, since $M$ is an interior p-submanifold of $\mathbb{R}_{1}^{N+1}$,
the vector fields $U,V$ can be extended to vector fields $\tilde{U},\tilde{V}$
on $\mathbb{R}_{1}^{N+1}$ tangent to the boundary. We then have $\overline{\nabla}_{U}V=\left(\nabla_{\tilde{U}}\tilde{V}\right)_{|M}^{\top}$.
Now, since $\tilde{U}$ and $\tilde{V}$ are tangent to the boundary
of $\mathbb{R}_{1}^{N+1}$, $\nabla_{\tilde{U}}\tilde{V}$ is tangent
to the boundary of $\mathbb{R}_{1}^{N+1}$ as well. In other words,
$\left\langle \nabla_{\tilde{U}}\tilde{V},\partial_{X}\right\rangle \equiv0$
along the boundary of $\mathbb{R}_{1}^{N+1}$. Since $M$ is orthogonal
to the boundary, it then follows that $\left\langle \overline{\nabla}_{U}V,\overline{\nabla}x\right\rangle _{\overline{g}}\equiv0$
along $\partial M$, which means exactly that $\overline{\nabla}_{U}V$
is tangent to the boundary of $\partial M$.
\end{proof}

\subsection{\label{subsec:Sufficiency}Sufficiency}

Here we prove the sufficiency part of Theorem \ref{thm:main-thm}.
More explicitly, we prove the following
\begin{thm}
\label{thm:sufficiency}Let $\left(M^{m},g\right)$ be an asymptotically
hyperbolic manifold, with $m\geq2$. Suppose that the symmetric $0$-$2$-tensor
$\Ric g+\left(m-1\right)g$ vanishes to second order along $\partial M$.
Then $\left(M,g\right)$ admits an isometric interior p-embedding
into $\HH^{\mathrm{N}_{m+1}+1}$, where $\mathrm{N}_{k}$ is the smallest
positive integer such that every closed Riemannian $k$-manifold can
be isometrically embedded in $\mathbb{R}^{\mathrm{N}_{k}}$.
\end{thm}

The proof builds upon that of the main theorem in \cite{UsulaIsometricEmbeddings},
whose statement is recalled in Theorem \ref{thm:previous-main-theorem}.
We now sketch the proof of this result, and build upon it later.

Suppose that $\left(M^{m},g\right)$ is conformally compact, with
sectional curvature at infinity satisfying $\kappa_{\infty}^{g}>-1$.
In order to isometrically p-embed $\left(M,g\right)$ into $\HH^{N+1}$
for some $N\geq m$, we need to find a boundary defining function
$x$ and a smooth map $\boldsymbol{v}:M\to\mathbb{R}^{N}$ such that
$\left(x,\boldsymbol{v}\right)$ is an embedding $M\to\mathbb{R}_{1}^{N+1}$,
and
\[
x^{2}g-dx^{2}=\left|d\boldsymbol{v}\right|^{2}.
\]
If we can construct $x$ so that $G=x^{2}g-dx^{2}$ is a smooth metric
on $M$, then the existence of $\boldsymbol{v}$ is guaranteed by
the Nash Embedding Theorem: it suffices to choose an isometric embedding
$\boldsymbol{v}:\left(M,G\right)\to\mathbb{R}^{N}$. To construct
such an $x$, choose a metric $h_{0}$ in the conformal infinity of
$g$, and let $r$ be a boundary defining function inducing the metric
$h_{0}$. Write $\overline{g}=r^{2}g$, define $K^{2}=\left|dr\right|_{\overline{g}}^{2}$,
and call $\overline{\nabla}r$ the gradient of $r$ with respect to
$\overline{g}$. In a sufficiently small collar $[0,\varepsilon)\times\partial M\to M$
of $\partial M$ in $M$ induced by the flow of $\overline{\nabla}r$,
we can write
\[
g=\frac{1}{K^{2}}\frac{dr^{2}}{r^{2}}+\frac{h\left(r\right)}{r^{2}}
\]
for some smooth family of metrics $h:[0,\varepsilon)\to C^{\infty}\left(\partial M;S^{2}\left(T^{*}\partial M\right)\right)$.
Note that in this case $K$ can be non-constant; however, since $\kappa_{\infty}^{g}=-K_{|\partial M}^{2}$,
we can choose $\varepsilon>0$ small enough so that $K^{2}<1$ in
the chosen collar. Now, choose a smooth function $Q:[0,\varepsilon)\to\left[0,1\right]$
such that $Q\left(0\right)=0$ and $Q\left(r\right)=1$ for $r\geq\varepsilon/2$.
Define
\[
\varphi\left(r\right)=-\int_{0}^{r}Q\left(t\right)\frac{dt}{t}.
\]
Define now $x=e^{\varphi}r$. By construction, $x$ is constant for
$r\geq\varepsilon/2$, and therefore it can be extended smoothly to
the whole of $M$; moreover, away from the smaller collar $r<\varepsilon/2$,
we have $G=x^{2}g$ and therefore $G$ is positive-definite there.
Finally, in the collar, we have
\[
G=e^{2\varphi}\left[\left(K^{-2}-\left(1-Q\right)^{2}\right)dr^{2}+h\left(r\right)\right]
\]
and the inequalities $K^{-2}>1$ and $\left(1-Q\right)^{2}\leq1$
guarantee that $G$ is positive-definite in the collar as well.

The same strategy does not work straight away when $g$ is asymptotically
hyperbolic. Indeed, in this case, the function $K$ equals $1$ along
$\partial M$, and for this reason the coefficient $e^{2\varphi}\left(K^{-2}-\left(1-Q\right)^{2}\right)$
of $dr^{2}$ in the formula above vanishes along $\partial M$. Nevertheless,
if $\Ric g+\left(m-1\right)g=O\left(x^{2}\right)$, the construction
above can still be used to produce an isometric interior p-embedding
$\left(x,\boldsymbol{v}\right):\left(M,g\right)\to\HH^{N+1}$ for
some $N$ sufficiently big. The key to the argument is the following
\begin{prop}
\label{prop:key-prop}Let $\left(M^{m},g\right)$ be an asymptotically
hyperbolic manifold. Suppose that $\Ric g+\left(m-1\right)g=O\left(x^{2}\right)$
as a section of $S^{2}\left(^{0}T^{*}M\right)$. Then there exists
a boundary defining function $x$, a compact Riemannian manifold with
boundary $\left(\mathcal{M}^{m+1},\mathcal{G}\right)$, and a smooth
map $F:M\to\mathcal{M}$ with the following properties:
\begin{enumerate}
\item the pull-back of a boundary defining function on $\mathcal{M}$ is
the square of a boundary defining function on $M$; in particular,
$F$ maps $M^{\circ}$ to $\mathcal{M}^{\circ}$ and $\partial M$
to $\partial\mathcal{M}$;
\item the restrictions $F^{\circ}:M^{\circ}\to\mathcal{M}^{\circ}$ and
$F_{\partial}:\partial M\to\partial\mathcal{M}$ are both embeddings;
\item $F^{*}\mathcal{G}=x^{2}g-dx^{2}$.
\end{enumerate}
\end{prop}

\begin{proof}
The argument starts with a variant of the previous construction. First
of all, since $g$ is asymptotically hyperbolic, given a metric $h_{0}\in\mathfrak{c}_{\infty}\left(g\right)$
we can choose a \emph{geodesic} boundary\emph{ }defining function
$r$ inducing $h_{0}$. Calling $\overline{g}=r^{2}g$, we can choose
a geodesic collar $[0,\varepsilon)\times\partial M\to M$ where $\left|dr\right|_{\overline{g}}^{2}\equiv1$
and
\[
g=\frac{dr^{2}+h\left(r\right)}{r^{2}}.
\]
Thanks to the decay condition $\Ric g+\left(m-1\right)g=O\left(r^{2}\right)$,
we have $h=h_{0}+O\left(r^{2}\right)$. Therefore, using Borel's Lemma,
we can decompose $h\left(r\right)=a\left(r^{2}\right)+r^{3}b\left(r^{2}\right)$,
where $a,b$ are two smooth functions $[0,\varepsilon^{2})\to C^{\infty}\left(\partial M;S^{2}\left(T^{*}\partial M\right)\right)$
with $a\left(0\right)=h_{0}$. This decomposition separates the even
terms in the expansion of $h$, contained in the expansion of $a\left(r^{2}\right)$,
from the odd terms, contained in the expansion of $r^{3}b\left(r^{2}\right)$.

We now choose a smooth function $q:[0,\varepsilon^{2})\to\left[0,1\right]$
with the following properties:
\begin{enumerate}
\item $q\left(0\right)=0$ and $q\left(s\right)>0$ for $s>0$;
\item $q'\left(0\right)>0$;
\item $q\left(s\right)\equiv1$ for $s\in[\left(\varepsilon/2\right)^{2},\varepsilon^{2})$.
\end{enumerate}
We define $Q\left(r\right)=q\left(r^{2}\right)$. As in the sketched
proof above, we define
\begin{align*}
\varphi\left(r\right) & :=-\int_{0}^{r}Q\left(t\right)\frac{dt}{t}\\
 & =-\int_{0}^{r^{2}}q\left(t\right)\frac{1}{2}\frac{dt}{t}\\
 & =:\psi\left(r^{2}\right).
\end{align*}
Observe that since $q\left(0\right)=0$, $\psi$ is a well-defined
smooth function $[0,\varepsilon^{2})\to\mathbb{R}$; therefore, $\varphi:[0,\varepsilon)\to\mathbb{R}$
is well-defined and smooth; moreover, $\varphi$ has even expansion
as $r\to0^{+}$ since $\varphi\left(r\right)=\psi\left(r^{2}\right)$.
As above, define $x=e^{\varphi}r$; by construction, $x$ is constant
on $r\geq\varepsilon/2$, and we can therefore extend it to the whole
of $M$. Moreover, since $q\left(0\right)=0$, we can write $q\left(s\right)=s\tilde{q}\left(s\right)$
for some smooth function $\tilde{q}:[0,\varepsilon^{2})\to\mathbb{R}$;
the condition $q'\left(0\right)>0$ then is equivalent to $\tilde{q}\left(0\right)>0$.

Call now $G=x^{2}g-dx^{2}$. Then, for $r\geq\varepsilon/2$ we have
$G=x^{2}g$, while in the chosen collar we have
\begin{align*}
G & =e^{2\varphi}\left[Q\left(2-Q\right)dr^{2}+h\left(r\right)\right]\\
 & =e^{2\psi\left(r^{2}\right)}\left[\tilde{q}\left(r^{2}\right)\left(2-r^{2}\tilde{q}\left(r^{2}\right)\right)r^{2}dr^{2}+a\left(r^{2}\right)+r^{3}b\left(r^{2}\right)\right].
\end{align*}
The idea is to express $G$ as the pull-back of a non-degenerate metric
on another compact manifold with boundary\emph{ $\mathcal{M}^{m+1}$}.
We first consider the manifold 
\[
U^{\left(1\right)}:=[0,\varepsilon^{2})\times\partial M\times\left(-\varepsilon^{3},\varepsilon^{3}\right),
\]
where $[0,\varepsilon^{2})$ is parametrized by the variable $s$
and $\left(-\varepsilon^{3},\varepsilon^{3}\right)$ is parametrized
by the variable $t$. Later, we shall see $U^{\left(1\right)}$ as
an open subset of a compact manifold with boundary $\mathcal{M}$.
We define on $U^{\left(1\right)}$ the symmetric $2$-tensor
\[
\mathcal{G}^{\left(1\right)}=e^{2\psi\left(s\right)}\left[\tilde{q}\left(s\right)\left(2-s\tilde{q}\left(s\right)\right)\frac{1}{4}ds^{2}+a\left(s\right)+tb\left(s\right)+\lambda^{2}dt^{2}-\frac{9}{4}s\lambda^{2}ds^{2}\right],
\]
where $\lambda^{2}\in\mathbb{R}^{+}$ is still to be determined. We
want to choose $\lambda^{2}$ in such a way that $\mathcal{G}^{\left(1\right)}$
is positive-definite. First of all, since $\left(a\left(s\right)+tb\left(s\right)\right)_{|s=t=0}=h_{0}$,
we can assume $\varepsilon$ small enough so that $a\left(s\right)+tb\left(s\right)$
is positive-definite for every $\left(s,t\right)\in[0,\varepsilon^{2})\times\left(-\varepsilon^{3},\varepsilon^{3}\right)$.
We can now just choose $\lambda^{2}$ so that the coefficient of $ds^{2}$
is strictly positive on $U^{\left(1\right)}$. This can be done, because
$q$ takes values in $\left[0,1\right]$ with $q\left(s\right)=0$
only for $s=0$, and $q\equiv1$ for $s\geq\left(\frac{\varepsilon}{2}\right)^{2}$.
These properties indeed guarantee that $\tilde{q}\left(2-s\tilde{q}\right)$
is bounded from below by a positive constant, and therefore we can
choose $\lambda^{2}$ so that
\[
\tilde{q}\left(2-s\tilde{q}\right)-9s\lambda^{2}>0
\]
for every $s\in[0,\varepsilon^{2})$. We now consider the map
\begin{align*}
F:[0,\varepsilon)\times\partial M & \to U^{\left(1\right)}\\
\left(r,p\right) & \mapsto\left(r^{2},p,r^{3}\right).
\end{align*}
Note that $s$ is a boundary defining function for $[0,\varepsilon^{2})\times\partial M\times\left(-\varepsilon^{3},\varepsilon^{3}\right)$,
and $F^{*}s=r^{2}$. Moreover, it is immediate to see that both restrictions
$F_{\partial}$ and $F^{\circ}$ are embeddings. Finally, elementary
computations show that $F^{*}\mathcal{G}^{\left(1\right)}=G$.

Now, as anticipated above, we extend $U^{\left(1\right)}=[0,\varepsilon^{2})\times\partial M\times\left(-\varepsilon^{3},\varepsilon^{3}\right)$
to a compact manifold with boundary $\mathcal{M}$; a way to do it
is to extend $[0,\varepsilon^{2})\times\partial M$ to $M$, and think
of $\left(-\varepsilon^{3},\varepsilon^{3}\right)$ as a neighborhood
of $1$ in $S^{1}$. We can as well extend $F$ to a smooth map $F:M\to\mathcal{M}$,
whose restriction to the interiors is an embedding $F^{\circ}:M^{\circ}\to\mathcal{M}^{\circ}$.
Indeed, it suffices to choose $F=\left(f,\tau\right)$, where $f:M\to M$
restricts to the map $[0,\varepsilon)\times\partial M\to[0,\varepsilon^{2})\times\partial M$
given by $\left(r,p\right)\mapsto\left(r^{2},p\right)$ and is a diffeomorphism
$M^{\circ}\to M^{\circ}$, and $\tau:M\to S^{1}$ is an arbitrary
smooth function extending the map $[0,\varepsilon)\to\left(-\varepsilon^{3},\varepsilon^{3}\right)$
given by $r\mapsto r^{3}$. By construction, the map $F$ satisfies
the first two points in the statement of the proposition. Moreover,
since $F^{\circ}:M^{\circ}\to\mathcal{M}^{\circ}$ is an embedding,
we can extend $G$ off $F^{\circ}\left(M^{\circ}\right)$ to a smooth
metric $\mathcal{G}^{\left(2\right)}$ on $U^{\left(2\right)}:=\mathcal{M}^{\circ}$
such that $F^{*}\mathcal{G}^{\left(2\right)}=G$. Finally, choose
an arbitrary metric $\mathcal{G}^{\left(3\right)}$ on $U^{\left(3\right)}:=\mathcal{M}\backslash F\left(M\right)$.
By means of a partition of unity subordinated to the open cover $\left\{ U^{\left(1\right)},U^{\left(2\right)},U^{\left(3\right)}\right\} $
of $\mathcal{M}$, we can then combine $\mathcal{G}^{\left(1\right)},\mathcal{G}^{\left(2\right)},\mathcal{G}^{\left(3\right)}$
into a metric $\mathcal{G}$ on $\mathcal{M}$ such that $F^{*}\mathcal{G}=G$.
\end{proof}
We can now easily conclude the proof of Theorem \ref{thm:sufficiency}:
\begin{proof}[Proof of Theorem \ref{thm:sufficiency}]
Choose $x$, $\mathcal{M}^{m+1}$, $\mathcal{G}$ and $F:M\to\mathcal{M}$
as in Proposition \ref{prop:key-prop}. Then we have
\[
F^{*}\mathcal{G}=x^{2}g-dx^{2}.
\]
Denote now by $2\mathcal{M}$ the smooth double of $\mathcal{M}$,
and extend $\mathcal{G}$ off $\mathcal{M}\subset2\mathcal{M}$ to
a smooth metric on $2\mathcal{M}$. Using the Nash Embedding Theorem,
we can find a smooth embedding $\boldsymbol{w}:2\mathcal{M}\to\mathbb{R}^{\mathrm{N}_{m+1}}$;
therefore, calling $\boldsymbol{v}$ the restriction of $\boldsymbol{w}$
to $\mathcal{M}$, we obtain $\boldsymbol{v}^{*}g_{\mathrm{eucl}}=\mathcal{G}$
or equivalently $\left|d\boldsymbol{v}\right|^{2}=\mathcal{G}$. It
follows that the map $\left(x,\boldsymbol{v}\circ F\right):M\to\mathbb{R}_{1}^{\mathrm{N}_{m+1}+1}$
pulls the hyperbolic $0$-metric 
\[
g_{\mathrm{hyp}}=\frac{dX^{2}+\left|d\boldsymbol{Y}\right|^{2}}{X^{2}}
\]
back to $g$. To conclude, we must show that $\left(x,\boldsymbol{v}\circ F\right)$
is an interior p-embedding. By construction, $\boldsymbol{v}\circ F:M\to\mathbb{R}^{\mathrm{N}_{m+1}}$
restricts to an embedding $M^{\circ}\to\mathbb{R}^{\mathrm{N}_{m+1}}$
and to an embedding $\partial M\to\mathbb{R}^{\mathrm{N}_{m+1}}$.
While $d\left(\boldsymbol{v}\circ F\right)$ is degenerate along $\partial M$
in the $\partial_{r}$ direction, the map $\left(x,\boldsymbol{v}\circ F\right)$
is non-degenerate at that direction because $dx\left(\partial_{r}\right)=1$
along $\partial M$. Therefore, $\left(x,\boldsymbol{v}\circ F\right)$
is an injective immersion, and hence an embedding by the compactness
of $M$. Finally, by construction, $\left(x,\boldsymbol{v}\circ F\right)$
pulls the boundary defining function $X$ for $\mathbb{R}_{1}^{\mathrm{N}_{m+1}+1}$
back to the boundary defining function $x$ for $M$; consequently,
$\left(x,\boldsymbol{v}\circ F\right)$ is an interior p-embedding,
concluding the proof.
\end{proof}
\bibliography{allpapers}{}
\bibliographystyle{alpha}
\newpage
\end{document}